# An introduction to the Barnes double gamma function with an application to an integral involving the cotangent function

Donal F. Connon

dconnon@btopenworld.com

17 February 2019

**Abstract**

We show a direct approach to an integral involving the cotangent function (which was originally discovered by Kinkelin in 1856 and published by him four years later in Crelle's Journal [14]). New derivations of Lerch's identity and the Gosper/Vardi functional equation are also presented.

## 1. Introduction

We have the Weierstrass canonical form of the gamma function [19, p.1]

(1.1) $$\frac{1}{\Gamma(x)} = xe^{\gamma x}\prod_{n=1}^{\infty}\left\{\left(1+\frac{x}{n}\right)e^{-\frac{x}{n}}\right\}$$

where $\gamma$ is Euler's constant defined by $\gamma = \lim_{n\to\infty}[H_n - \log n]$ and $H_n$ are the harmonic numbers $H_n = \sum_{k=1}^{n}\frac{1}{k}$.

Taking logarithms results in

(1.2) $$\log\Gamma(1+x) = -\gamma x - \sum_{n=1}^{\infty}\left[\log\left(1+\frac{x}{n}\right)-\frac{x}{n}\right]$$

and letting $x=0$ in (1.2) we immediately see that $\log\Gamma(1)=0$, and hence $\Gamma(1)=1$.

The Barnes double gamma function $G(x)$ may be defined by [19, p.25]

(1.3) $$\log G(1+x) = \frac{1}{2}x\log(2\pi) - \frac{1}{2}x(1+x) - \frac{1}{2}\gamma x^2 + \sum_{n=1}^{\infty}\left[\frac{1}{2n}x^2 - x + n\log\left(1+\frac{x}{n}\right)\right]$$

and letting $x=0$ we immediately see that $\log G(1)=0$, and hence we have $G(1)=1$.

The above definitions of $\Gamma(x)$ and $G(x)$ will be employed throughout this paper and we shall only assume knowledge of those properties of the functions which may be deduced from the respective definitions.

As shown below, using these representations, it is relatively straightforward to derive Kinkelin's integral [14] which is valid for $-1<u<1$

(1.4) $$\int_0^u \pi x \cot \pi x \, dx = u \log(2\pi) + \log \frac{G(1-u)}{G(1+u)}$$

## 2. Some basic properties of the gamma function $\Gamma(x)$

We first look at some of the basic properties of the gamma function.

We see from (1.2) that

(2.1) $$\log \Gamma(2) = -\gamma - \sum_{n=1}^{\infty} \left[ \log\left(1+\frac{1}{n}\right) - \frac{1}{n} \right]$$

Euler's constant $\gamma$ may be defined as the limit of the sequence

$$\gamma = \lim_{n\to\infty} [H_n - \log n]$$

$$= \lim_{n\to\infty} [H_n - \log(n+1) + \log(n+1) - \log n]$$

$$= \lim_{n\to\infty} [H_n - \log(n+1)] + \lim_{n\to\infty} \log\left(1+\frac{1}{n}\right)$$

$$= \lim_{n\to\infty} [H_n - \log(n+1)]$$

We have

$$\sum_{k=1}^{n} \log\left(1+\frac{1}{k}\right) = \sum_{k=1}^{n} [\log(k+1) - \log k]$$

and it is easily seen that this telescopes to

$$\sum_{k=1}^{n} \log\left(1+\frac{1}{k}\right) = \log(n+1)$$

We can therefore write Euler's constant $\gamma$ as the infinite series

(2.2) $$\gamma = \sum_{n=1}^{\infty} \left[ \frac{1}{n} - \log\left(1+\frac{1}{n}\right) \right]$$

and we immediately see from (2.1) that $\log \Gamma(2) = 0$ and hence we deduce that $\Gamma(2) = 1$.

□

The combination of (1.2) with (2.2) results in

$$\log \Gamma(1+x) = \sum_{n=1}^{\infty} \left[ x \log\left(1+\frac{1}{n}\right) - \log\left(1+\frac{x}{n}\right) \right]$$

and differentiation gives us (an infrequently used formula)

$$\psi(1+x) = \sum_{n=1}^{\infty} \left[ \log\left(1+\frac{1}{n}\right) - \frac{1}{n+x} \right]$$

which is equivalent to the well-known representation [19, p.14]

$$\psi(x) = \lim_{N\to\infty} \left[ \log N - \sum_{n=0}^{N} \frac{1}{n+x} \right]$$

□

We see from (1.2) that

(2.3) $$\log \Gamma(1+x) + \log \Gamma(1-x) = -\sum_{n=1}^{\infty} \left[ \log\left(1+\frac{x}{n}\right) + \log\left(1-\frac{x}{n}\right) \right]$$

$$= -\sum_{n=1}^{\infty} \log\left(1-\frac{x^2}{n^2}\right)$$

and assuming the well-known identity originally postulated by Euler

(2.4) $$\sin \pi x = \pi x \prod_{n=1}^{\infty} \left(1-\frac{x^2}{n^2}\right)$$

we then obtain Euler's reflection formula for the gamma function [19, p.3]

(2.5) $$\Gamma(x)\Gamma(1-x) = \frac{\pi}{\sin \pi x}$$

Logarithmic differentiation of (2.5) results in

(2.6) $$\psi(x) - \psi(1-x) = -\pi \cot \pi x$$

where $\psi(x)$ is the digamma function defined by $\psi(x) := \frac{d}{dx} \log \Gamma(x)$.

Using the expansion for the digamma function (obtained by differentiating (1.2))

(2.7) $$\psi(1+x) = -\gamma + \sum_{n=1}^{\infty} \left( \frac{1}{n} - \frac{1}{n+x} \right)$$

we have upon letting $1+x \to x$

$$\psi(x) = -\gamma + \sum_{n=1}^{\infty}\left(\frac{1}{n} - \frac{1}{n-1+x}\right)$$

$$= -\gamma + \sum_{m=0}^{\infty}\left(\frac{1}{m+1} - \frac{1}{m+x}\right)$$

$$= -\gamma + 1 - \frac{1}{x} + \sum_{m=1}^{\infty}\left(\frac{1}{m+1} - \frac{1}{m+x}\right)$$

$$= -\gamma + 1 - \frac{1}{x} + \sum_{m=1}^{\infty}\left(\frac{1}{m} - \frac{1}{m+x} + \frac{1}{m+1} - \frac{1}{m}\right)$$

This gives us

(2.8) $$\psi(x) = -\gamma - \frac{1}{x} + \sum_{m=1}^{\infty}\left(\frac{1}{m} - \frac{1}{m+x}\right)$$

and we therefore obtain the functional equation

(2.9) $$\psi(1+x) = \psi(x) + \frac{1}{x}$$

Integrating this over the interval $[1,u]$ results in

$$\log\Gamma(1+u) - \log\Gamma(2) = \log\Gamma(u) - \log\Gamma(1) + \log u$$

and, since $\Gamma(2) = \Gamma(1) = 1$, we deduce the familiar functional equation

(2.10) $$\Gamma(1+u) = u\,\Gamma(u)$$

**3. Raabe's integral**

We consider the integral

(3.1) $$u(x) = \int_x^{x+1} \log\Gamma(t)\,dt$$

and note the derivative

$$u'(x) = \log\Gamma(x+1) - \log\Gamma(x)$$

$$= \log x$$

where we have used (2.10). Hence, upon integration we get

$$u(x) = x\log x - x + a$$

where $a$ is the constant of integration and we see that

$$u(0)=\int_0^1 \log\Gamma(t)\,dt = a$$

Letting $t \to 1-t$ we see that

$$\int_0^1 \log\Gamma(t)\,dt = \int_0^1 \log\Gamma(1-t)\,dt = a$$

and therefore we have

$$\int_0^1 \log\left[\Gamma(t)\Gamma(1-t)\right]dt = 2a$$

Then, using Euler's reflection formula (2.5) $\Gamma(t)\Gamma(1-t)=\dfrac{\pi}{\sin \pi t}$ , we get

$$\int_0^1 \log\left[\frac{\pi}{\sin \pi t}\right]dt = 2a$$

and thus we find

$$a=\frac{1}{2}\log\pi - \frac{1}{2}\int_0^1 \log\sin \pi t\,dt$$

We have the well-known integral which was first evaluated by Euler

$$\int_0^{\pi/2} \log\sin x\,dx = -\frac{\pi}{2}\log 2$$

and hence we have

$$\int_0^{\pi} \log\sin x\,dx = \int_0^{\pi/2} \log\sin x\,dx + \int_{\pi/2}^{\pi} \log\sin x\,dx = -\pi\log 2\,.$$

We therefore obtain the integration constant $a=\frac{1}{2}\log(2\pi)$ and accordingly we have for $x>0$

(3.2) $$u(x)=\int_x^{x+1} \log\Gamma(t)\,dt = x\log x - x + \frac{1}{2}\log(2\pi)$$

which also applies in the limit as $x \to 0$ to give

(3.3) $$u(0)=\int_0^1 \log\Gamma(t)\,dt = \frac{1}{2}\log(2\pi)$$

These identities were first shown by Joseph Ludwig Raabe (1801-1859) in Crelle's Journal in 1840 (see [18] and [22, p.261]). The formula (3.2) can be used to derive Stirling's asymptotic formula for $\Gamma(x)$ and the factorial $n!$: for example, using the first mean-value theorem for integrals we have

(3.4) $$\int_x^{x+1} \log\Gamma(t)\,dt = \log\Gamma(\xi) = x\log x - x + \frac{1}{2}\log(2\pi)$$

where $x < \xi < x+1$. As a first approximation, one could take $\xi = x$ which would give us the asymptotic formula

$$\log\Gamma(x) \approx x\log x - x + \frac{1}{2}\log(2\pi)$$

In fact, the even more truncated approximation $\log\Gamma(x) \approx x\log x$ is good enough to be successfully employed in various applications of theoretical physics [17, p.167 & 219].

**4. Some basic properties of the double gamma function $G(x)$**

Differentiating (1.3) gives us

$$\frac{G'(1+x)}{G(1+x)} = \frac{1}{2}\log(2\pi) - \frac{1}{2} - x - \gamma x + \sum_{n=1}^{\infty}\left[\frac{1}{n}x - 1 + \frac{n}{n+x}\right]$$

which may be written as

(4.1) $$\frac{G'(1+x)}{G(1+x)} = \frac{1}{2}\log(2\pi) - \frac{1}{2} - x - \gamma x + x\sum_{n=1}^{\infty}\left[\frac{1}{n} - \frac{1}{n+x}\right]$$

Using (2.7) we obtain

(4.2) $$\frac{G'(1+x)}{G(1+x)} = \frac{1}{2}\log(2\pi) - \frac{1}{2} - x + x\psi(1+x)$$

We now integrate this over the interval $[0,u]$ to obtain (noting that $\log G(1) = 0$)

(4.2.1) $$\log G(1+u) = \frac{1}{2}u\log(2\pi) - \frac{1}{2}u(u+1) + \int_0^u x\psi(1+x)\,dx$$

and integration by parts gives us

$$\int_0^u x\psi(1+x)\,dx = u\log\Gamma(1+u) - \int_0^u \log\Gamma(1+x)\,dx$$

Hence, we obtain Alexeiewsky's theorem [19, p.32]

(4.3) $$\int_0^u \log\Gamma(1+x)\,dx = \frac{1}{2}[\log(2\pi)-1]u - \frac{1}{2}u^2 + u\log\Gamma(1+u) - \log G(1+u)$$

Alternatively, we may also compute $\int_0^u \log\Gamma(1+x)\,dx$ directly using the representation (1.2) as follows:

$$\int_0^u \log\Gamma(1+x)\,dx = -\frac{1}{2}\gamma u^2 - \sum_{n=1}^{\infty}\left[(n+u)\log\left(1+\frac{u}{n}\right) - u - \frac{u^2}{2n}\right]$$

$$= -\frac{1}{2}\gamma u^2 - \sum_{n=1}^{\infty}\left[n\log\left(1+\frac{u}{n}\right) - u + \frac{u^2}{2n}\right] - u\sum_{n=1}^{\infty}\left[\log\left(1+\frac{u}{n}\right) - \frac{u}{n}\right]$$

and using (1.2) and (1.3) we see that (4.3) follows directly.

With $u=1$ in (4.3) we have

$$\log G(2) = \frac{1}{2}\log(2\pi) - 1 - \int_0^1 \log\Gamma(1+x)\,dx$$

$$= \frac{1}{2}\log(2\pi) - 1 - \int_0^1 \log\Gamma(x)\,dx - \int_0^1 \log x\,dx$$

and, referring to Raabe's integral (3.3), we deduce that $\log G(2) = 0$ and thus we see that $G(2) = 1$.

□

**Proposition:**

(4.4) $$G(1+x) = \Gamma(x)G(x).$$

**Proof:**

Letting $x \to 1+x$ in (4.2) gives us

$$\frac{G'(2+x)}{G(2+x)} = \frac{1}{2}\log(2\pi) - \frac{3}{2} - x + (1+x)\psi(2+x)$$

$$= \frac{1}{2}\log(2\pi) - \frac{3}{2} - x + x\psi(2+x) + \psi(2+x)$$

$$= \frac{1}{2}\log(2\pi) - \frac{3}{2} - x + x\psi(1+x) + \frac{x}{1+x} + \psi(1+x) + \frac{1}{1+x}$$

$$= \frac{1}{2}\log(2\pi) - \frac{1}{2} - x + x\psi(1+x) + \psi(1+x)$$

$$= \frac{G'(1+x)}{G(1+x)} + \psi(1+x)$$

Thus we have

$$\frac{G'(2+x)}{G(2+x)} - \frac{G'(1+x)}{G(1+x)} = \psi(1+x)$$

and integration gives us

$$\log G(2+x) - \log G(2) - \log G(1+x) = \log \Gamma(1+x)$$

Letting $1+x \to x$, and noting that $\log G(2) = 0$, we see that

$$\log G(1+x) - \log G(x) = \log \Gamma(x)$$

## 5. An integral of the cotangent function

Referring to (4.2)

$$\frac{G'(1+x)}{G(1+x)} = \frac{1}{2}\log(2\pi) - \frac{1}{2} - x + x\psi(1+x)$$

and with $x \to -x$ we have

$$\frac{G'(1-x)}{G(1-x)} = \frac{1}{2}\log(2\pi) - \frac{1}{2} + x - x\psi(1-x)$$

We then have

$$\frac{G'(1+x)}{G(1+x)} + \frac{G'(1-x)}{G(1-x)} = \log(2\pi) - 1 + x[\psi(1+x) - \psi(1-x)]$$

From (2.6) and (2.9) we see that

$$\psi(1+x) - \psi(1-x) = \frac{1}{x} - \pi \cot \pi x$$

and hence we have

(5.1) $$\frac{G'(1+x)}{G(1+x)} + \frac{G'(1-x)}{G(1-x)} = \log(2\pi) - \pi x \cot \pi x$$

Integration of this results in Kinkelin's integral

(5.2) $$\int_0^u \pi x \cot \pi x \, dx = u \log(2\pi) + \log \frac{G(1-u)}{G(1+u)}$$

Using (4.4), $G(1+u) = G(u)\Gamma(u)$, we may immediately deduce Euler's integral

$$\int_0^{1/2} \pi x \cot \pi x \, dx = -\int_0^{1/2} \log \sin \pi x \, dx = \frac{1}{2}\log 2$$

or equivalently

$$\int_0^{\pi/2} \log \sin x \, dx = -\frac{\pi}{2}\log 2$$

□

We also see from (4.1) that

$$\frac{G'(1+x)}{G(1+x)} + \frac{G'(1-x)}{G(1-x)} = \log(2\pi) - 1 + x\sum_{n=1}^{\infty}\left[\frac{1}{n-x} - \frac{1}{n+x}\right]$$

$$= \log(2\pi) - 1 + 2x^2\sum_{n=1}^{\infty}\frac{1}{n^2 - x^2}$$

and, using (5.1), we deduce the decomposition formula for the cotangent function

$$\pi x \cot \pi x = 1 + 2x^2 \sum_{n=1}^{\infty} \frac{1}{x^2 - n^2}$$

Integrating this gives us

$$\int_0^u \left[\pi \cot \pi x - \frac{1}{x}\right] dx = \int_0^u \sum_{n=1}^{\infty} \frac{2x}{x^2 - n^2} \, dx$$

The Weierstrass $M$ test shows that we may integrate the series term by term to obtain

$$\log \frac{\sin \pi u}{\pi u} = \sum_{n=1}^{\infty} \log\left(1 - \frac{u^2}{n^2}\right)$$

and hence we obtain (2.4)

$$\sin \pi u = \pi u \prod_{n=1}^{\infty}\left(1 - \frac{u^2}{n^2}\right)$$

## 6. An elementary derivation of Lerch's identity

Lerch [16] established the following relationship between the gamma function and the Hurwitz zeta function in 1894

(6.1) $$\varsigma'(0, x) = \log \Gamma(x) - \frac{1}{2}\log(2\pi)$$

where $\varsigma(s,x)$ is the Hurwitz zeta function defined initially for $\operatorname{Re} s > 1$ and $x > 0$ by

(6.2) $$\varsigma(s,x) = \sum_{n=0}^{\infty} \frac{1}{(n+x)^s}$$

and $\varsigma'(s,x) := \frac{\partial}{\partial s}\varsigma(s,x)$. Note that $\varsigma(s,x)$ may be analytically continued to the whole $s$-plane except for a simple pole at $s=1$. For example, Hasse (1898-1979) showed that [13]

(6.3) $$\varsigma(s,x) = \frac{1}{s-1}\sum_{n=0}^{\infty} \frac{1}{n+1}\sum_{k=0}^{n} \binom{n}{k} \frac{(-1)^k}{(k+x)^{s-1}}$$

is a globally convergent series for $\varsigma(s,x)$ and, except for $s=1$, provides an analytic continuation of $\varsigma(s,x)$ to the entire complex plane.

It may be noted that $\varsigma(s,1) = \varsigma(s)$.

**Proof:**

Since the Hurwitz zeta function $\varsigma(s,x)$ is analytic in the whole complex plane except for $s=1$, its partial derivatives commute in the region where the function is analytic: we therefore have

$$\frac{\partial}{\partial x}\frac{\partial}{\partial s}\varsigma(s,x) = \frac{\partial}{\partial s}\frac{\partial}{\partial x}\varsigma(s,x)$$

$$= -\frac{\partial}{\partial s}[s\varsigma(s+1,x)]$$

Hence we have

(6.4) $$\frac{\partial}{\partial x}\frac{\partial}{\partial s}\varsigma(s,x) = -\varsigma(s+1,x) - s\frac{\partial}{\partial s}\varsigma(s+1,x)$$

Differentiating again with respect to $x$ results in

(6.5) $$\frac{\partial^2}{\partial x^2}\frac{\partial}{\partial s}\varsigma(s,x) = (2s+1)\varsigma(s+2,x) + s(s+1)\varsigma'(s+2,x)$$

and with $s=0$ in (6.5) we have

$$\frac{d^2}{dx^2}\varsigma'(0,x) = \varsigma(2,x)$$

Reference to (2.8) and (6.2) shows us that $\varsigma(2,x) = \psi'(x)$ and hence we have

$$\frac{d^2}{dx^2}[\varsigma'(0,x)-\log\Gamma(x)]=0$$

Integration results in

$$\varsigma'(0,x)-\log\Gamma(x)=ax+b$$

where $a$ and $b$ are integration constants to be obtained. With $x=1$ we have $\varsigma'(0)=a+b$ and $x=2$ gives us $\varsigma'(0,2)=2a+b$.

From the series definition (6.2) of the Hurwitz zeta function it is easily seen that

$$\varsigma(s,x)=\varsigma(s,1+x)+\frac{1}{x^s}$$

and hence we have upon differentiation with respect to $s$

$$\varsigma'(s,x)=\varsigma'(s,1+x)-\frac{\log x}{x^s}$$

Thus we see that $\varsigma'(s,1)=\varsigma'(s,2)$ and therefore $\varsigma'(0,2)=\varsigma'(0)$. Hence we determine that $a=0$. Accordingly, we obtain

(6.6) $$\varsigma'(0,x)-\log\Gamma(x)=b$$

A little more effort is required to determine the integration constant $b$. We use the known identity

(6.7) $$\varsigma\left(s,\frac{1}{2}\right)=\left[2^s-1\right]\varsigma\left(s\right)$$

and differentiation gives us

(6.8) $$\varsigma'\left(s,\frac{1}{2}\right)=\left[2^s-1\right]\varsigma'\left(s\right)+2^s\log 2\varsigma\left(s\right)$$

Letting $s=0$ shows that

$$\varsigma'\left(0,\frac{1}{2}\right)=\varsigma\left(0\right)\log 2$$

We see from (6.10) below that $\varsigma(0,1)=-\frac{1}{2}$ and thus $\varsigma\left(0\right)=-\frac{1}{2}$. Hence, we have

$$\varsigma'\left(0,\frac{1}{2}\right)=-\frac{1}{2}\log 2$$

Then using the well-known result [19, p.3] $\Gamma\left(\frac{1}{2}\right) = \sqrt{\pi}$, which may be obtained by letting $x = \frac{1}{2}$ in (2.5), we determine by letting $x = \frac{1}{2}$ in (6.6) that $b = -\frac{1}{2}\log(2\pi)$. Therefore, we have Lerch's identity

$$\varsigma'(0,x) = \log\Gamma(x) - \frac{1}{2}\log(2\pi)$$

Other fairly elementary derivations of Lerch's identity are given in [6] and [8]; more convoluted ones are given in [19, p.91] and [22, p.271].

With $x = 1$ we see that

(6.9) $$\varsigma'(0) = -\frac{1}{2}\log(2\pi)$$ □

Letting $s = 1 - m$ in (6.3) gives us

$$\varsigma(1-m,x) = -\frac{1}{m}\sum_{n=0}^{\infty}\frac{1}{n+1}\sum_{k=0}^{n}\binom{n}{k}(-1)^k (k+x)^m$$

and, following in the footsteps of F. Lee Cook [15], we showed in [10] that the Bernoulli polynomials could be represented by

$$B_m(x) = \sum_{n=0}^{\infty}\frac{1}{n+1}\sum_{k=0}^{n}\binom{n}{k}(-1)^k (k+x)^m$$

Hence we obtain another derivation of the well-known result [5, p.264]

(6.10) $$\varsigma(1-m,x) = -\frac{B_m(x)}{m}$$

Since $\sum_{k=0}^{n}\binom{n}{k}(-1)^k (k+x)^m = 0$ for $m > n = 0,1,2,...$ we therefore have the finite polynomial expression

(6.11) $$B_m(x) = \sum_{n=0}^{m}\frac{1}{n+1}\sum_{k=0}^{n}\binom{n}{k}(-1)^k (k+x)^m$$

□

Differentiating (6.3) gives us

(6.12) $$(s-1)\varsigma'(s,x) + \varsigma(s,x) = -\sum_{n=0}^{\infty}\frac{1}{n+1}\sum_{k=0}^{n}\binom{n}{k}(-1)^k \frac{\log(k+x)}{(k+x)^{s-1}}$$

and with $s=0$ we have

$$-\varsigma'(0,x)+\varsigma(0,x)=-\sum_{n=0}^{\infty}\frac{1}{n+1}\sum_{k=0}^{n}\binom{n}{k}(-1)^k(k+x)\log(k+x)$$

Employing Lerch's identity (6.1) results in

(6.13) $$\log\Gamma(x)=\sum_{n=0}^{\infty}\frac{1}{n+1}\sum_{k=0}^{n}\binom{n}{k}(-1)^k(k+x)\log(k+x)+\frac{1}{2}-x+\frac{1}{2}\log(2\pi)$$

where we have used $\varsigma(0,x)=-B_1(x)$.

□

Letting $s=-1$ in (6.4) gives us

$$\frac{\partial}{\partial x}\varsigma'(-1,x)=-\varsigma(0,x)+\varsigma'(0,x)$$

$$=B_1(x)+\log\Gamma(x)-\frac{1}{2}\log(2\pi)$$

Integration over the interval $[1,u]$ gives us

$$\varsigma'(-1,u)-\varsigma'(-1,1)=u\log\Gamma(u)-\log G(1+u)$$

where we have used (3.3) and (4.3). This is equivalent to

(6.14) $$G(1+u)-u\log\Gamma(u)=\varsigma'(-1)-\varsigma'(-1,u)$$

We shall see a slightly different derivation in (7.1) below.

□

Letting $s=-1$ in (6.12) gives us

$$-2\varsigma'(-1,x)+\varsigma(-1,x)=-\sum_{n=0}^{\infty}\frac{1}{n+1}\sum_{k=0}^{n}\binom{n}{k}(-1)^k(k+x)^2\log(k+x)$$

Using (6.14) and $\varsigma(-1,x)=-\frac{B_2(x)}{2}$ we obtain

(6.15)

$$\log G(1+x)=-\frac{1}{2}\sum_{n=0}^{\infty}\frac{1}{n+1}\sum_{k=0}^{n}\binom{n}{k}(-1)^k(k+x)^2\log(k+x)+x\log\Gamma(x)+\frac{1}{4}B_2(x)+\varsigma'(-1)$$

The above expansions for $\log\Gamma(x)$ and $\log G(1+x)$ were previously obtained in [8] in a rather more circuitous manner.

Equation (6.15) may be written as

$$\log G(1+x) - x\log\Gamma(x) = -\frac{1}{2}\sum_{n=0}^{\infty}\frac{1}{n+1}\sum_{k=0}^{n}\binom{n}{k}(-1)^k (k+x)^2 \log(k+x) + \frac{1}{4}B_2(x) + \varsigma'(-1)$$

and we note that this corresponds with (7.1) below.

## 7. The Gosper/Vardi functional equation

**Proposition 7.1**

(7.1) $$G(1+x) - x\log\Gamma(x) = \varsigma'(-1) - \varsigma'(-1,x)$$

The functional equation (7.1) was derived by Vardi [20] in 1988 and also by Gosper [12] in 1997.

**Proof:**

With $s=-1$ in (6.4) we have

$$\frac{d}{dx}\varsigma'(-1,x) = -\varsigma(0,x) + \varsigma'(0,x)$$

Then, using (6.10) we see that

$$\varsigma(0,x) = \frac{1}{2} - x$$

This identity may also be obtained directly from (6.3). We then have

$$\frac{d}{dx}\varsigma'(-1,x) = x - \frac{1}{2} + \log\Gamma(x) - \frac{1}{2}\log 2\pi$$

and integration over the interval $[1,x]$ results in

$$\varsigma'(-1,x) - \varsigma'(-1) = \frac{1}{2}x(x-1) + \int_1^x \log\Gamma(x)\,dx - \frac{1}{2}(x-1)\log 2\pi$$

since $\varsigma'(-1,1) = \varsigma'(-1)$. Hence, using Alexeiewsky's theorem (4.3), we obtain

$$\log G(1+x) - x\log\Gamma(x) = \varsigma'(-1) - \varsigma'(-1,x)$$

We can use this result to determine the value of $G\left(\frac{1}{2}\right)$. Letting $s=-1$ in (6.8) gives us

$$\varsigma'\left(-1,\frac{1}{2}\right) = -\frac{1}{2}\varsigma'(-1) + \frac{1}{2}\log 2\varsigma(-1)$$

and we have $\varsigma\left(-1\right) = -\frac{B_2(1)}{2} = -\frac{B_2}{2} = -\frac{1}{12}$ .

We therefore find the known result

(7.2) $$\log G\left(\frac{1}{2}\right) = \frac{1}{24}\log 2 - \frac{1}{4}\log\pi + \frac{3}{2}\varsigma'\left(-1\right)$$

With $s = -2$ in (6.5) we have

$$\frac{\partial^2}{\partial x^2}\varsigma'(-2,x) = -3\varsigma(0,x) + 2\varsigma'(0,x)$$

and using (4.3) we may deduce the integral of $\log G(1+x)$ in terms of $\varsigma'(-2,x)$ .

**Proposition 7.2**

(7.3) $$n\int_0^x \varsigma'(1-n,u)\,du = \frac{B_{n+1} - B_{n+1}(x)}{n(n+1)} + \varsigma'(-n,x) - \varsigma'(-n)$$

This integral for $n \geq 1$ was obtained by Adamchik [2] in 1998.

**Proof:**

We recall (6.4)

$$\frac{\partial}{\partial u}\frac{\partial}{\partial s}\varsigma(s,u) = -\varsigma(s+1,u) - s\frac{\partial}{\partial s}\varsigma(s+1,u)$$

and, upon integrating this over $[1,x]$ we see that

$$-s\int_1^x \varsigma'(s+1,u)du = \int_1^x \frac{\partial}{\partial u}\frac{\partial}{\partial s}\varsigma(s,u)du + \int_1^x \varsigma(s+1,u)\,du$$

We therefore get

(7.4) $$-s\int_1^x \varsigma'(s+1,u)du = \varsigma'(s,x) - \varsigma'(s,1) + \int_1^x \varsigma(s+1,u)\,du$$

It should be noted that we were reluctant to integrate over the interval $[0,x]$ because we would then end up with $\varsigma'(-n,0)$ which does not appear to be defined.

With $s = -n$ we have

(7.5) $$n\int_1^x \varsigma'(1-n,u)du = \varsigma'(-n,x) - \varsigma'(-n) + \int_1^x \varsigma(1-n,u)\,du$$

Then using $\varsigma(1-n,u) = -\dfrac{B_n(u)}{n}$ for $n \ge 1$ we obtain

(7.6) $$n\int_1^x \varsigma'(1-n,u)\,du = \frac{B_{n+1} - B_{n+1}(x)}{n(n+1)} + \varsigma'(-n,x) - \varsigma'(-n)$$

and with $x = 2$ we have

$$n\int_1^2 \varsigma'(1-n,u)\,du = \frac{B_{n+1} - B_{n+1}(2)}{n(n+1)} + \varsigma'(-n,2) - \varsigma'(-n)$$

It is well known [19] that

$$B_{n+1}(1+x) = B_{n+1}(x) + (n+1)x^n$$

and we therefore have $B_{n+1}(2) = B_{n+1} + (n+1)$. Thus we obtain

(7.7) $$n\int_1^2 \varsigma'(1-n,u)\,du = -\frac{1}{n} + \varsigma'(-n,2) - \varsigma'(-n)$$

Differentiating

$$\varsigma(s,1+x) = \varsigma(s,x) - \frac{1}{x^s}$$

we obtain

$$\varsigma'(s,1+x) = \varsigma'(s,x) + \frac{\log x}{x^s}$$

and, in particular, we have

(7.8) $$\varsigma'(1-n,1+x) = \varsigma'(1-n,x) + x^{n-1}\log x$$

With $x = 1$ we see that

$$\varsigma'(1-n,2) = \varsigma'(1-n)$$

We now integrate (7.8) over the interval $[0,1]$

(7.9) $$\int_0^1 [\varsigma'(1-n,1+x) - \varsigma'(1-n,x)]\,dx = \int_0^1 x^{n-1}\log x\,dx$$

and, as shown below, parametric differentiation is of assistance. We see that

$$\int_0^u x^p \log x\, dx = \frac{\partial}{\partial p}\int_0^u x^p\, dx$$

$$= \frac{\partial}{\partial p}\frac{u^{p+1}}{p+1}$$

$$= \frac{u^{p+1}[(p+1)\log u - 1}{(p+1)^2}$$

Therefore, we obtain the known integral

$$\int_0^1 x^{n-1} \log x\, dx = -\frac{1}{n^2}$$

Hence we have

(7.10) $$\int_0^1 [\varsigma'(1-n,1+x) - \varsigma'(1-n,x)]\, dx = -\frac{1}{n^2}$$

An obvious change of variables gives us

$$\int_0^1 \varsigma'(1-n,1+x)\, dx = \int_1^2 \varsigma'(1-n,t)\, dt$$

$$\int_0^1 \varsigma'(1-n,1+x)\, dx = \frac{1}{n}\left[-\frac{1}{n} + \varsigma'(-n,2) - \varsigma'(-n)\right]$$

and (7.10) and (7.7) gives us

$$\frac{1}{n}\left[-\frac{1}{n} + \varsigma'(-n,2) - \varsigma'(-n)\right] - \int_0^1 \varsigma'(1-n,u)\, du = -\frac{1}{n^2}$$

Since $\varsigma'(-n,2) = \varsigma'(-n)$ we immediately see that

(7.11) $$\int_0^1 \varsigma'(1-n,u)\, du = 0$$

This enables us to write (7.6) as

$$n\int_0^x \varsigma'(1-n,u)\, du = \frac{B_{n+1} - B_{n+1}(x)}{n(n+1)} + \varsigma'(-n,x) - \varsigma'(-n)$$

This proves the proposition. Adamchik [2] stated that this integral is valid for $n \geq 1$ but, as shown below, we see that it is valid for $n = 0$.

With $n = 0$ we have

(7.12) $$\int_0^x \varsigma'(0,u)\,du = \frac{1}{2}[B_2 - B_2(x)] + \varsigma'(-1,x) - \varsigma'(-1)$$

$$= \frac{1}{2}(x - x^2) + \varsigma'(-1,x) - \varsigma'(-1)$$

$$= \frac{1}{2}(x - x^2) + x \log\Gamma(x) - \log G(1+x)$$

and this leads to another derivation of Alexeiewsky's theorem (4.3).

**Proposition 7.3**

(7.13) $$\int_0^x \varsigma''(0,u)du = \varsigma''(-1,x) - \varsigma''(-1) + x - x^2 + 2[x\log\Gamma(x) - \log G(1+x)]$$

**Proof**

We differentiate (6.4) to obtain

$$\frac{\partial}{\partial u}\frac{\partial^2}{\partial s^2}\varsigma(s,u) = -2\frac{\partial}{\partial s}\varsigma(s+1,u) - s\frac{\partial^2}{\partial s^2}\varsigma(s+1,u)$$

and, upon integrating this over $[1,x]$ we see that

$$-s\int_1^x \varsigma''(s+1,u)du = \int_1^x \frac{\partial}{\partial u}\frac{\partial^2}{\partial s^2}\varsigma(s,u)du + 2\int_1^x \frac{\partial}{\partial s}\varsigma(s+1,u)\,du$$

We therefore get

(7.14) $$-s\int_1^x \varsigma''(s+1,u)du = \varsigma''(s,x) - \varsigma''(s) + 2\int_1^x \frac{\partial}{\partial s}\varsigma(s+1,u)\,du$$

With $s = -1$ we have

$$\int_1^x \varsigma''(0,u)du = \varsigma''(-1,x) - \varsigma''(-1) + 2\int_1^x \varsigma'(0,u)\,du$$

Employing (7.12) we obtain

$$\int_1^x \varsigma''(0,u)du = \varsigma''(-1,x) - \varsigma''(-1) + x - x^2 + 2[x\log\Gamma(x) - \log G(1+x)]$$

We note from differentiating

$$\varsigma(s,1+x) = \varsigma(s,x) - \frac{1}{x^s}$$

that

$$\lim_{x\to 0+}\left[\varsigma^{(n)}(-m,x) - (-1)^n x^m \log^n x\right] = \varsigma^{(n)}(-m)$$

Using mathematical induction together with L'Hôpital's rule, it is easy to prove that $\lim_{x\to 0+}\left(x\log^n x\right) = 0.$

Hence, provided $m \geq 1$ we have

$$\text{(7.15)} \qquad \lim_{x\to 0+} \varsigma^{(n)}(-m,x) = \varsigma^{(n)}(-m)$$

This explains why $\lim_{x\to 0}\varsigma'(0,x)$ does not exist whereas $\lim_{x\to 0}\varsigma'(-1,x)$ does exist. Therefore, we have

$$\int_0^x \varsigma''(0,u)du = \varsigma''(-1,x) - \varsigma''(-1) + x - x^2 + 2[x\log\Gamma(x) - \log G(1+x)]$$

and we see that

$$\int_0^1 \varsigma''(0,u)du = 0$$

We showed in [24] that

$$\int_0^x \varsigma''(0,u)\,du - \varsigma''(0)x = x[\log^2 x - 2\log x + 2]$$

$$+\sum_{n=1}^{\infty}\left[(n+x)\left\{\log^2(n+x) - 2\log(n+x) + 2\right\} - n\left\{\log^2 n - 2\log n + 2\right\} - x\log^2 n - \frac{1}{2}x^2[\log^2(n+1) - \log^2 n]\right]$$

and we obtain

$$\varsigma''(-1,x) - \varsigma''(-1) + x - x^2 + 2[x\log\Gamma(x) - \log G(1+x)] = \varsigma''(0)x + x[\log^2 x - 2\log x + 2]$$

$$+\sum_{n=1}^{\infty}\left[(n+x)\left\{\log^2(n+x) - 2\log(n+x) + 2\right\} - n\left\{\log^2 n - 2\log n + 2\right\} - x\log^2 n - \frac{1}{2}x^2[\log^2(n+1) - \log^2 n]\right]$$

With $x = 1$ we obtain

$$\varsigma''(0) = -2 - \sum_{n=1}^{\infty} \left[ (n+1)\log^2(n+1) - n\log^2 n - 2[(n+1)\log(n+1) - n\log n] + 2 - \frac{1}{2}[\log^2(n+1) + \log^2 n] \right]$$

*WolframAlpha* confirms that the above series is convergent.

Using the finite telescoping sum

$$\sum_{k=1}^{N} [a_{k+1} - a_k] = a_{N+1} - a_1$$

we have

$$\sum_{n=1}^{N} \left[ (n+1)\log^2(n+1) - n\log^2 n - 2[(n+1)\log(n+1) - n\log n] + 2 - \frac{1}{2}[\log^2(n+1) + \log^2 n] \right]$$

$$= (N+1)\log^2(N+1) - 2(N+1)\log(N+1) + 2N - \frac{1}{2}\log^2(N+1)$$

$$= \left( N + \frac{1}{2} \right) \log^2(N+1) - 2(N+1)\log(N+1) + 2N$$

We therefore obtain

$$\varsigma''(0) = -2 + \lim_{N\to\infty} \left[ -\left( N + \frac{1}{2} \right) \log^2(N+1) + 2(N+1)\log(N+1) - 2N \right]$$

and comparing this with [24]

$$\varsigma''(0) = \lim_{N\to\infty} \left[ \sum_{k=1}^{N} \log^2 k - \left( N + \frac{1}{2} \right) \log^2 N + 2N\log N - 2N \right]$$

we obtain

$$2 = \lim_{N\to\infty} \left[ -\left( N + \frac{1}{2} \right) [\log^2(N+1) - \log^2 N] + 2(N+1)\log(N+1) - 2N\log N - \sum_{k=1}^{N} \log^2 k \right]$$

**Proposition 7.3**

(7.13) $$n\int_1^x \varsigma''(1-n,u)\,du = \varsigma''(-n,x) - \varsigma''(-n) + \frac{2}{n}\left[ \frac{B_{n+1} - B_{n+1}(x)}{n(n+1)} + \varsigma'(-n,x) - \varsigma'(-n) \right]$$

**Proof**

With $s = -n$ in (7.14) we have

$$n\int_1^x \varsigma''(1-n,u)du = \varsigma''(-n,x) - \varsigma''(-n) + 2\int_1^x \varsigma'(1-n,u)\,du$$

and substituting (7.3) gives us

$$n\int_1^x \varsigma''(1-n,u)du = \varsigma''(-n,x) - \varsigma''(-n) + \frac{2}{n}\left[\frac{B_{n+1} - B_{n+1}(x)}{n(n+1)} + \varsigma'(-n,x) - \varsigma'(-n)\right]$$

We see that

$$\int_0^1 \varsigma''(1-n,u)du = 0$$

**8. Some connections with the Clausen function** $\mathrm{Cl}_2(x)$

In [9] we proved the basic identity

(8.1) $$\int_a^b p(x)\cot(\alpha x/2)\,dx = 2\sum_{n=1}^{\infty}\int_a^b p(x)\sin\alpha nx\,dx$$

where we initially required that $p(x)$ is a twice continuously differentiable function. Equation (8.1) is valid provided (i) $\sin(\alpha x/2) \neq 0 \ \forall \ x \in [a,b]$ or, alternatively, (ii) if $\sin(\alpha\eta/2) = 0$ for some $\eta \in [a,b]$ then $p(\eta) = 0$.

As later shown in [11], (8.1) is actually valid for a wider class of suitably behaved functions (which are not necessarily twice continuously differentiable).

With $p(x) = x$ and $\alpha = 2\pi$ in (8.1) we have for $-1 < u < 1$

$$\int_0^u \pi x\cot\pi x\,dx = 2\pi\sum_{n=1}^{\infty}\int_0^u x\sin 2\pi nx\,dx$$

and since

$$\int_0^u x\sin 2\pi nx\,dx = \frac{\sin 2\pi nu}{(2\pi n)^2} - \frac{u\cos 2\pi nu}{2\pi n}$$

we obtain

(8.2) $$\int_0^u \pi x\cot\pi x\,dx = -u\sum_{n=1}^{\infty}\frac{\cos 2n\pi u}{n} + \frac{1}{2\pi}\sum_{n=1}^{\infty}\frac{\sin 2n\pi u}{n^2}$$

Equating this with (5.2) we obtain

(8.3) $$u\log(2\pi) + \log\frac{G(1-u)}{G(1+u)} = -u\sum_{n=1}^{\infty}\frac{\cos 2n\pi u}{n} + \frac{1}{2\pi}\sum_{n=1}^{\infty}\frac{\sin 2n\pi u}{n^2}$$

Using the familiar trigonometric series shown in Carslaw's book [7, p.241]

$$\log\left[2\sin(\pi u)\right] = -\sum_{n=1}^{\infty}\frac{\cos 2n\pi u}{n} \qquad (0 < u < 1)$$

we obtain

$$u\log(2\pi) + \log\frac{G(1-u)}{G(1+u)} = u\log\left[2\sin(\pi u)\right] + \frac{1}{2\pi}\sum_{n=1}^{\infty}\frac{\sin 2n\pi u}{n^2}$$

which, using (2.5), may be written as

$$u\log[\Gamma(u)\Gamma(1-u)] + \log\frac{G(1-u)}{G(1+u)} = \frac{1}{2\pi}\sum_{n=1}^{\infty}\frac{\sin 2n\pi u}{n^2}$$

We recall the Clausen function $\mathrm{Cl}_2(u)$ defined by

$$\mathrm{Cl}_2(u) = \sum_{n=1}^{\infty}\frac{\sin nu}{n^2} = -\int_0^u \log\left(2\sin\frac{\theta}{2}\right)d\theta$$

and we see that

(8.4) $$\frac{1}{2\pi}\mathrm{Cl}_2(2\pi u) = u\log[\Gamma(u)\Gamma(1-u)] - \log\frac{G(1+u)}{G(1-u)}$$

as previously shown by Adamchik [3].

From the series definition of the Clausen function $\mathrm{Cl}_2(u)$ we easily see that

$$\frac{1}{2}\mathrm{Cl}_2(2\pi u) = \mathrm{Cl}_2(\pi u) - \mathrm{Cl}_2(\pi(1-u))$$

and (8.4) gives us

$$\frac{1}{\pi}[\mathrm{Cl}_2(\pi u) - \mathrm{Cl}_2(\pi(1-u))] = u\log[\Gamma(u)\Gamma(1-u)] - \log\frac{G(1+u)}{G(1-u)}$$

$$= u\log[\Gamma(u)\Gamma(1-u)] + \log G(1-u) - \log G(u) - \log\Gamma(u)$$

$$= u\log\Gamma(1-u) - (1-u)\log\Gamma(u) + \log G(1-u) - \log G(u)$$

Hence we see that

(8.5)
$$\frac{1}{\pi}\mathrm{Cl}_2(\pi u) + (1-u)\log\Gamma(u) + \log G(u) = \frac{1}{\pi}\mathrm{Cl}_2(\pi(1-u)) + u\log\Gamma(1-u) + \log G(1-u)$$

Defining $f(u)$ as

$$f(u)=\frac{1}{\pi}\text{Cl}_2(\pi u)+(1-u)\log\Gamma(u)+\log G(u)$$

we see that $f(u)=f(1-u)$.

Using (7.1) we see that

$$G(1+u)-u\log\Gamma(u)=\varsigma'(-1)-\varsigma'(-1,u)$$

or equivalently

$$G(u)+(1-u)\log\Gamma(u)=\varsigma'(-1)-\varsigma'(-1,u)$$

Hence (8.5) gives us

$$\frac{1}{\pi}\text{Cl}_2(\pi u)-\varsigma'(-1,u)=\frac{1}{\pi}\text{Cl}_2(\pi(1-u))-\varsigma'(-1,1-u)$$

and employing $\frac{1}{2}\text{Cl}_2(2\pi u)=\text{Cl}_2(\pi u)-\text{Cl}_2(\pi(1-u))$ we obtain the known result [1]

$$\frac{1}{2\pi}\text{Cl}_2(2\pi u)=\varsigma'(-1,u)-\varsigma'(-1,1-u)$$

It was also shown in [9] that

$$\int_a^b\frac{p(x)}{\sin\alpha x}dx=2\sum_{n=0}^{\infty}\int_a^b p(x)\sin\left[(2n+1)\alpha x\right]dx \tag{8.6}$$

Equation (8.1) is valid provided (i) $\sin(\alpha x)\neq 0\ \forall\ x\in[a,b]$ or, alternatively, (ii) if $\sin(\alpha\eta)=0$ for some $\eta\in[a,b]$ then $p(\eta)=0$.

Since

$$\int_0^u x\sin(2n+1)\pi x\,dx=-\frac{u\cos(2n+1)u}{(2n+1)\pi}+\frac{\sin(2n+1)u}{(2n+1)^2\pi^2}$$

we obtain

$$\int_0^u\frac{\pi x}{\sin\pi x}dx=2\pi\sum_{n=0}^{\infty}\left[-\frac{u\cos(2n+1)u}{(2n+1)\pi}+\frac{\sin(2n+1)u}{(2n+1)^2\pi^2}\right]$$

Wang [21] has shown that (see also [19, p.32])

$$\int_0^u\frac{\pi x}{\sin\pi x}dx=u\log(2\pi)+\log\frac{G(1+u)}{G(1-u)}-4\log\frac{G(1+u/2)}{G(1-u/2)} \tag{8.7}$$

and we therefore have

$$u\log(2\pi)+\log\frac{G(1+u)}{G(1-u)}-4\log\frac{G(1+u/2)}{G(1-u/2)}=2\pi\sum_{n=0}^{\infty}\left[-\frac{u\cos(2n+1)u}{(2n+1)\pi}+\frac{\sin(2n+1)u}{(2n+1)^2\pi^2}\right]$$

**9. Open access to our own work**

This paper contains references to various other papers and, rather surprisingly, most of them are currently freely available on the internet. Surely now is the time that all of our work should be freely accessible by all. The mathematics community should lead the way on this by publishing everything on arXiv, or in an equivalent open access repository. We think it, we write it, so why hide it? You know it makes sense.

**REFERENCES**

[1] V.S. Adamchik, A Class of Logarithmic Integrals. Proceedings of the 1997 International Symposium on Symbolic and Algebraic Computation. ACM, Academic Press, 1-8, 2001. http://www-2.cs.cmu.edu/~adamchik/articles/logs.htm

[2] V.S. Adamchik, Polygamma Functions of Negative Order. J. Comp. and Applied Math.100, 191-199, 1998. http://www.cs.cmu.edu/~adamchik/articles/polyg.pdf

[3] V.S. Adamchik, On the Barnes Function. Proceedings of the 2001 International Symposium on Symbolic and Algebraic Computation. (July 22-25, 2001, London, Canada), Academic Press, 15-20, 2001 http://www-2.cs.cmu.edu/~adamchik/articles/issac01/issac01.pdf

[4] V.S. Adamchik, Contributions to the Theory of the Barnes Function. Computer Physics Communications, 2003. http://www-2.cs.cmu.edu/~adamchik/articles/barnes1.pdf

[5] T.M. Apostol, Introduction to Analytic Number Theory. Springer-Verlag, New York, Heidelberg and Berlin, 1976.

[6] B.C. Berndt, The Gamma Function and the Hurwitz Zeta Function. Amer. Math. Monthly, 92,126-130, 1985.

[7] H.S. Carslaw, Introduction to the theory of Fourier Series and Integrals. Third Ed. Dover Publications Inc, 1930.

[8] D.F. Connon, Some series and integrals involving the Riemann zeta function, binomial coefficients and the harmonic numbers. Volume II(a), 2007. arXiv:0710.4023 [pdf]

[9] D.F. Connon, Some series and integrals involving the Riemann zeta function, binomial coefficients and the harmonic numbers. Volume V, 2007. arXiv:0710.4047 [pdf]

[10] D.F. Connon, Some series and integrals involving the Riemann zeta function, binomial coefficients and the harmonic numbers. Volume VI, 2007. arXiv:0710.4032 [pdf]

[11] D.F. Connon, Some applications of the Dirichlet integrals to the summation of series and the evaluation of integrals involving the Riemann zeta function. 2012. arXiv:1212.0441 [pdf]

[12] R.W. Gosper, $\int_{n/4}^{m/6} \log \Gamma(z)dz$. In Special Functions, q-series and related topics. Amer. Math. Soc., Vol. 14.

[13] H. Hasse, Ein Summierungsverfahren für Die Riemannsche $\varsigma$ - Reithe. Math.Z.32, 458-464, 1930.
http://dz-srv1.sub.uni-goettingen.de/sub/digbib/loader?ht=VIEW&did=D23956&p=462

[14] V.H. Kinkelin, Üeber eine mit der Gammafunction verwandte Transcendente und deren Anwendung auf die Integralrechnung.
Journal für die reine und angewandte Mathematik, 1860, 57, 122-158.
http://www.digizeitschriften.de/no_cache/home/open-access/nach-zeitschriftentiteln/

[15] F. Lee Cook, A Simple Explicit Formula for the Bernoulli Numbers. The Two-Year College Mathematics Journal 13, 273-274, 1982.

[16] M. Lerch, Dalsi studie v oboru Malmstenovskych rad, Rozpravy Ceske Akad. 3, no. 28, 1894, 63 pp.

[17] M.S. Longair, Theoretical Concepts in Physics: An Alternative View of Theoretical Reasoning in Physics, Cambridge University Press, 1984.

[18] J.L. Raabe, Angenäherte Bestimmung der Factorenfolge

$1.2.3.4.5...n = \Gamma(1+n) = \int x^n e^{-x} dx$, wenn $n$ eine sehr gross Zahl ist.

Journal für die reine und angewandte Mathematik (25),146 -159, 1843.
http://www.digizeitschriften.de/en/dms/toc/?PPN=PPN243919689_0025

[19] H.M. Srivastava and J. Choi, Series Associated with the Zeta and Related Functions. Kluwer Academic Publishers, Dordrecht, the Netherlands, 2001.

[20] I. Vardi, Determinants of Laplacians and multiple gamma functions, SIAM J. Math. Anal. 19 (1988) 493–507.

[21] X. Wang, The Barnes $G$ -function and the Catalan constant. Kyushu J. Math. 67 (2013), 105-116

[22] E.T. Whittaker and G.N. Watson, A Course of Modern Analysis: An Introduction to the General Theory of Infinite Processes and of Analytic Functions; With an Account of the Principal Transcendental Functions. Fourth Ed., Cambridge University Press, Cambridge, London and New York, 1963.

[23] N.-Y. Zhang and K.S. Williams, Some results on the generalized Stieltjes constants. Analysis, 14, 147-162 (1994). www.math.carleton.ca/~williams/papers/pdf/187.pdf

[24] D.F. Connon, Some new formulae involving the Stieltjes constants, 2019. https://arxiv.org/abs/1902.00510

Wessex House,
Devizes Road,
Upavon,
Pewsey,
Wiltshire SN9 6DL